\documentclass{amsart}
\usepackage{amsfonts}
\usepackage{amssymb}
\usepackage{times}
\font\got=eufb10
\def\m{\frac{1}{2}}
\def\k{\hbox{K}_{10}}
\def\aut{\mathop{\hbox{\rm Aut}}}
\def\GL{\mathop{\hbox{\rm GL}}}
\def\SL{\mathop{\hbox{\rm SL}}}
\def\T{\mathop{\mathcal T}}
\def\N{\mathop{\mathcal N}}

\def\WW{\hbox{\got W}}
\def\Z{\mathbb Z}
\def\span#1{\langle #1\rangle}
\def\ot{\otimes}
\def\l{\lambda}
\def\m{\mu}
\def\s{\sigma}
\def\ee{e\otimes e}
\def\xy{x\otimes y}
\def\yx{y\otimes x}
\def\xx{x\otimes x}
\def\yy{y\otimes y}
\def\xe{x\otimes e}
\def\ex{e\otimes x}
\def\ye{y\otimes e}
\def\ey{e\otimes y}

\newtheorem{pr}{Proposition}
\newtheorem{co}{Corollary}
\newtheorem{de}{Definition}
\newtheorem{te}{Theorem}
\begin{document}

\title{Gradings on the Kac superalgebra}

\thanks{ The   authors are
 partially supported by the MCYT grant MTM2007-60333 and by the Junta de
Andaluc\'{\i}a grants FQM-336, FQM-1215. The first author is also
supported by the PCI of the UCA `Teor\'\i a de Lie y Teor\'\i a de
Espacios de Banach' and by the PAI with project number
 FQM-298.}

\author[A. J. Calder\'on]{Antonio Jes\'us Calder\'on Mart\'{\i}n}
\email{ajesus.calderon@uca.es}
\address{Dpto. Matem\'aticas\\Facultad de Ciencias, Universidad de C\'adiz\\
Campus de Puerto Real, 11510, C\'adiz, Spain}

\author[C. Draper]{Cristina Draper Fontanals}
\email{cdf@uma.es}
\address{Dpto. Matem\'atica Aplicada\\Escuela Superior de Ingenier\'{\i}a Industrial, Universidad de M\'alaga\\
Campus de El Ejido, 29013, M\'alaga, Spain}

\author[C. Mart\'{\i}n]{C\'andido Mart\'{\i}n Gonz\'alez}
\email{candido@apncs.cie.uma.es}
\address{Dpto. \'Algebra, Geometr\'\i a y Topolog\'\i a\\Facultad de Ciencias, Universidad de M\'alaga\\
   Campus de Teatinos, 29080, M\'alaga, Spain}

\begin{abstract}
We describe the group gradings on the $K_{10}$ Jordan
superalgebra. There are 21 nonequivalent gradings, two of them
fine and 6 nontoral gradings.

\end{abstract}

\maketitle

\section{Introduction}

As a motivation for Lie superalgebras, the notion of supersymmetry
in theoretical physics reflects the known symmetry between bosons
and fermions. The mathematical structure formalizing this idea is
that of supergroup, or $\Z_2$-graded Lie group. As mentioned in
\cite{raptis}, their job is that of modelling  continuous
supersymmetry transformations between bosons and fermions. As Lie
algebras consist of generators of Lie groups, the infinitesimal
Lie group elements tangent to the identity, so $\Z_2$-graded Lie
algebras, otherwise known as Lie superalgebras, consist of
generators of (or infinitesimal) supersymmetry transformations.
Closely related to these (for instance by the Kantor-Koecher-Tits
construction) are the Jordan superalgebras (see \cite{consuelo}).
We are primarily interested in $\k$.

The debut of the simple exceptional Jordan superalgebra $\k$ in
the mathematical literature was in V. Kac's classification of
finite dimensional simple Jordan superalgebras over fields of
characteristic zero \cite{kac}. Though one can introduce $\k$ by
giving its multiplication table relative to some basis, more
conceptual approaches are possible. One of   such convenient
viewpoints is the given in \cite{benkart}.  Works as recent as
\cite{Daniel} and \cite{splittest} deal with the right definition
of the algebra over arbitrary rings of scalars (agreeing with the
usual $\k$ when the base ring is an algebraically closed field of
characteristic not $2$). The Grassmann envelope of the algebra has
been studied in \cite{grassmann}.  The interest on group gradings
of superalgebras seems to start with the work \cite{tval}.  The
paper \cite{Elduque} is essential for our study since it describes
the automorphism group of $\k$. Its relevance stems from the fact
that, in our setting, gradings are just the simultaneous
diagonalization  relative to commuting sets of automorphisms.

\section{Preliminaries}

Allthrough this work the base field $F$ will be an algebraically
closed field of zero characteristic. Let $J=J_0\oplus J_1$ be a
superalgebra over $F$. The term grading will always mean group
grading, that is, a decomposition in vector subspaces
$J=\oplus_{g\in G}J^g$ where $G$ is a finitely generated
 abelian group and the
 homogeneous spaces verify $J^gJ^h\subset J^{gh}$ (denoting by juxtaposition
  the product in $G$). We assume also that $G$
 is generated by the set of all $g$ such that $J^g\ne 0$, usually called the \emph{support} of the grading, and that
 the grading is compatible with the grading
$J=J_0\oplus J_1$ of the superalgebra $J$. This means that any
homogeneous component $J^i$ splits as $J^i=J_0^i\oplus J_1^i$
where $J_k^i:=J^i\cap J_k$ for $k\in\{0,1\}$. To  distinguish the
$\Z_2$-grading providing the superalgebra structure of $J$ from
the rest of its possible gradings, we denote it with subscripts
rather than with supscripts.

 Given
two gradings $J=\oplus_{g\in G}U^g$ and $J=\oplus_{h\in H}V^h$ we
shall say that they are {\em isomorphic} if there is a group
isomorphism $\sigma\colon G\to H$ and a (superalgebra)
automorphism $\varphi\colon J\to J$ such that
$\varphi(U^g)=V^{\sigma(g)}$ for all $g\in G$. We recall that a
superalgebra automorphism $\varphi$ is just an automorphism such
that the even and the odd part of the superalgebra are
$\varphi$-invariant. The above two gradings are said to be {\em
equivalent} if there are: (1) a bijection $\sigma\colon I\to I'$
between the supports of the first and second gradings
respectively, and, (2) a superalgebra automorphism $\varphi$ of
$J$ such that $\varphi(U^g)=V^{\sigma(g)}$ for any $g\in I$.

 Consider the $10$-dimensional $F$-algebra $\k$
whose basis is $(e,v_1,v_2,v_3,v_4,f,x_1,x_2,$ $y_1,y_2)$   with
multiplication table:
\smallskip

\begin{center}
\begin{tabular}{|c|cccccc||cccc|}
\hline
$\cdot$ &  $e$ & $v_1$ & $v_2$ & $v_3$ & $v_4$ & $f$ & $x_1$ & $x_2$ & $y_1$ & $y_2$ \cr
\hline
 $e$ &  $e$ & $v_1$ & $v_2$ & $v_3$ & $v_4$ & $0$ & $\frac{1}{2} x_1$ & $\frac{1}{2} x_2$ &
$\frac{1}{2} y_1$ & $\frac{1}{2} y_2$ \cr
 $v_1$ &  $v_1$ & $0$ & $2e$ & $0$ & $0$ & $0$  & $0$ & $0$ & $x_2$ & $-x_1$\cr
 $v_2$ &  $v_2$ & $2e$ & $0$ & $0$ & $0$ & $0$ & $-y_2$  & $y_1$ & $0$ & $0$ \cr $v_3$ & $v_3 $ & $0$ & $0$ & $0$ & $2e$ & $0$ & $0$ & $x_1$ & $y_2$ & $0$\cr
 $v_4$ & $v_4$ & $0$ & $0$ & $2e$ & $0$ & $0$ & $x_2$ & $0$ & $0$ & $y_1$\cr
 $f$ & $0$ & $0$ & $0$ & $0$ & $0$ & $f$ & $\frac{1}{2} x_1$ & $\frac{1}{2} x_2$ &
$\frac{1}{2} y_1$ & $\frac{1}{2} y_2$ \cr \hline $x_1$ &
$\frac{1}{2} x_1$ & $0$ & $-y_2$ & $0$ & $x_2$ & $\frac{1}{2} x_1$
& $0$ & $v_1$ & $e-3f$ & $v_3$\cr $x_2$ & $\frac{1}{2} x_2$ & $0$
& $y_1$ & $x_1$ & $0$ & $\frac{1}{2} x_2$ & $-v_1$ & $0$ & $v_4$ &
$e-3f$\cr $y_1$ & $\frac{1}{2} y_1$ & $x_2$ & $0$ & $y_2$ & $0$ &
$-\frac{1}{2} y_1$ & $3f-e$ & $-v_4$ & $0$ & $v_2$\cr $y_2$ &
$\frac{1}{2} y_2$ & $-x_1$ & $0$ & $0$ & $y_1$ & $\frac{1}{2} y_2$
& $-v_3$ & $3f-e$ & $-v_2$ & $0$\cr \hline
\end{tabular}
\end{center}
\smallskip

This is a Jordan superalgebra called the   {\em Kac superalgebra},
with even part generated by $(e,v_1,\ldots,v_4,f)$ and odd part
generated by $(x_1,x_2,y_1,y_2)$. We shall denote it by $\k$ as
usual in the mathematical literature.  We have used the basis
introduced in \cite{Racine} for this superalgebra, however a more
conceptual approach is possible (see for instance \cite{Elduque}.)
In order to adhere to this alternative approach we shall need the
{\em Kaplansky superalgebra}.  This is the 3-dimensional Jordan
superalgebra $K = K_0\oplus K_1$, with $K_0 = Fe$ and $K_1Fx\oplus Fy$, and with multiplication given by $e^2 = e$, $ex
=\frac{1}{2}x = xe$, $ey =\frac{1}{2}y = ye$, $xy = e =-y x$,
$x^2=y^2=0$. Following \cite{Elduque} we define on $K$ the
following supersymmetric bilinear form: $(e|e) = \frac{1}{2}$,
$(x|y) = 1$, $(K_0,K_1) = 0$, (it must be understood that the form
is symmetric on $K_0$ and alternating on $K_1$.) Consider now the
$F$-vector space over $F\cdot 1 \oplus (K \otimes_F K)$ and define
on it the product where $1$ is the identity element and
\begin{equation}\label{one}
(a\otimes b)(c\otimes d) = (-1)^{\bar b\bar c} (a c\otimes b d
-\frac{3}{4} (a,c)(b,d)\cdot 1)
\end{equation}
 for $a, b, c, d\in K$
homogeneous elements, where for any homogeneous element $x\in
K_0\cup K_1$, $\bar x=0$ if $x\in K_0$ and $\bar x=1$ if $x\in
K_1$. Then, a result of G. Benkart and A. Elduque states that $\k$
is isomorphic to $F\cdot 1\oplus(K \otimes K)$ with the above
super-product, by means of the isomorphism given in \cite[Theorem
2.1, p. 4]{Elduque}. Unless otherwise stated, we shall identify
$\k$ with the algebra whose product is (\ref{one}). Under this
identification the even part is $\span{1, e\ot e, x\ot x, x\ot y,
y\ot x, y\ot y}$ and the odd one is $\span{e\ot x, x\ot e, e\ot y,
y\ot e}$.

 In this work we shall have
the ocassion to consider automorphism groups of several Jordan
superalgebras. These are, of course, linear algebraic groups
whence some  aspects of this theory will be used in the sequel. We
must take into account that for a Jordan superalgebra $J=J_0\oplus
J_1$ the notation $\aut(J)$ will mean the group of all
  automorphisms (as above preserving the
homogeneous components.) The gradings on the Kaplansky
superalgebra $K$ may be computed rather straightforwardly from
scratch but we prefer to develop some algebraic group tools which
can be applied not only to $K$ but to other algebras (alternative,
Lie, Jordan, super-Jordan) of more complex nature. At this point
there are two key aspects to mention:
\begin{itemize}
\item Any grading is induced by a finitely generated abelian subgroup of
diagonalizable automorphisms of the automorphism group of the
algebra under study. The homogeneous components are the simultaneous
eigenspaces relative to the given group of automorphisms.
\item Any such subgroup is contained in the normalizer of some
maximal torus of the automorphism group of the algebra. This is an
algebraic group version of the Borel-Serre theorem for Lie groups.
It has been given by V. P. Platonov in Theorem~6 and Theorem~3.15,
p.~92 of \cite{Platonov}: A supersoluble subgroup of semisimple
elements of an algebraic group $G$ is contained in the normalizer
of a maximal torus. Here we must recall that a group is called
supersolvable (or supersoluble) if it has an invariant normal
series whose factors are all cyclic. Any finitely generated
abelian group is supersolvable. Since we are considering gradings
on Jordan superalgebras over this class of abelian groups, we may
assume that the set of automorphisms inducing the grading (as
simultaneous diagonalization) is contained in the normalizer of a
maximal torus of its automorphism group.
\end{itemize}

A special kind of gradings arises when we consider the inducing
automorphisms not only in the normalizer of a maximal torus, but
in the torus itself.

\begin{de}\rm
A grading of a superalgebra is  said to be   {\em toral} if it is
produced by   automorphisms within  a torus of the automorphism
group of the superalgebra.
\end{de}

Returning to the Kaplansky superalgebra $K$, it is straightforward
that any element in $\aut(K)$ fixes $e$ so that $\aut(K)$ can be
identified with a subgroup of $\GL_2(F)$.  Moreover taking into
account that $f(x)f(y)=e$, we easily check that $\aut(K) \cong
\SL_2(F)$. We are denoting by $\T$ the maximal torus of $\aut(K)$
(identified once and for all with $\SL_2(F)$) consisting of all
its diagonal matrices. This maximal torus is isomorphic to
$F^\times$ and a generic element in $\T$ will be denoted by
$t_\lambda:=\hbox{diag}(\lambda, \lambda^{-1})$ with $\lambda\in
F^\times$. The normalizer of $\T$ in $\SL_2(F)$ is the quasitorus
$\N=\T\cup\T\sigma$ where 
$\sigma= {\tiny 
\begin{pmatrix}
0 & 1\cr -1 &
0\end{pmatrix}}$ which corresponds to the automorphism of $K$ given by
$x\mapsto y\mapsto -x$  (of order $4$.) Thus, it is easy to see
that $\N/\T\cong\Z_2$ (since $\sigma^2=t_{-1}\equiv -1\in\T$.)
\begin{pr}
Any abelian subgroup of $\N$ is toral.
\end{pr}
Proof. It is essential the fact that $\sigma t=t^{-1}\sigma$ for
any $t\in\T$. Let $S$ be a nontrivial abelian subgroup of $\N$. It
$S\subset\T$ we are done. On the contrary, since $S\not\subset\T
\sigma$, we must have toral elements $t\in\T\cap S$ and also
elements $t'\sigma\in\T\sigma\cap S$. Since they must commute we
can write $t t' \sigma=t'\sigma t=t't^{-1}\sigma$ from which we
get $t^2=1$ and hence $t=\pm 1$. So $S\cap\T\subset\{\pm1\}$. On
the other hand, if $S$ contains two elements $t_1\sigma$ and
$t_2\sigma$, from the fact that they commute one gets $t_1=\pm
t_2$. Thus $S\subset\{\pm 1,\pm t\sigma\}$ for some $t\in\T$. But
  for any $t\in\T$ there is some
$t_1\in \T$ such that $t_1\sigma t t_1^{-1}=\sigma$ (take
$t_1^2=t$.) So by conjugating $S$ with $t_1$ we have
$t_1St_1^{-1}\subset\{\pm 1,\pm\sigma\}$, but this set is toral
since $\sigma$ is toral: for some $p\in\SL_2(F)$ we have $p\sigma
p^{-1}={\tiny\begin{matrix}i & 0\cr 0 &-i\end{matrix}}=t_i$ (so $p\{\pm
1,\pm\sigma\}p^{-1}\subset\T$.)
\begin{co}
All the gradings on the Kaplansky superalgebra $ K$ are toral. Up
to equivalence the nontrivial ones are the following:
\begin{itemize}
\item The $\Z_2$-grading $K=K_0\oplus K_1$ providing its superalgebra
structure.
\item The fine $\Z$-grading $K=K^{-1}\oplus K^0\oplus K^1$ such
that $K^{-1}=F x$, $K^0=F e$ and $K^1=F y$.
\end{itemize}
\end{co}
Proof. Consider a minimal set $S\ne\{1\}$ of diagonalizable
automorphisms inducing the grading. Since any two maximal tori of
the automorphism group are conjugated, we may assume $S\subset\T$.
Now, if $S$ contains some $t_\lambda$ with $\lambda\ne\pm 1$, we
get the fine $\Z$-grading. On the contrary, $S=\span{t_{-1}}$
induces the $\Z_2$-grading.

\section{Gradings on $\k$}

Now we develop a similar program for the $\k$ superalgebra. Since
our study depends heavily on the knowledge of $\aut(\k)$, we must
return to the reference \cite{Elduque} which gives full details on
this group. Firstly, if we take $f,g\in\aut(K)$ then we may define
an automorphism of $\k$ (denoted $(f,g)$) such that $(f,g)\colon
x\otimes y\mapsto f(x)\otimes g(y)$. Thus we have a group
monomorphism of $\aut(K)^2$ to $\aut(\k)$.
 But as stated in
 \cite{Elduque} there is an automorphism $\delta$ of $\k$ such that
 $\delta(x\otimes y)=(-1)^{\bar x\bar y}y\otimes x$. Moreover
 $\aut(\k)\cong\aut(K)^2\rtimes \{1,\delta\}\cong
\SL_2(F)^2\cup\SL_2(F)^2\delta$ where the product in
 $\SL_2(F)^2$ (or in $\aut(K)^2$) is componentwise and $(f,g)\delta=\delta(g,f)$. So
 $\aut(\k)$ has two connected components and the component
 of the unit is $\aut(\k)_0=\SL_2(F)^2$. A maximal torus of
 $\aut(\k)$ is then $\T^2=\T\times\T$ and its normalizer
 in $\aut(\k)$ is $\N(\T^2)=\N^2\cup\N^2\delta$, where we recall that
$\N$ is the normalizer of $\T$ in $\SL_2(F)$. From Section II we
know that the set of diagonalizable automorphisms of $\k$
producing any grading is (up to conjugacy) contained in
$\N(\T^2)$.

\subsection{Toral gradings on $\k$}

Any toral grading on $\k$  is
 isomorphic to a grading produced by a subquasitorus of $\T^2$.
We are denoting $t_{\lambda,\mu}:=(t_\lambda, t_\mu)\in \T^2$ for
any $\lambda,\mu\in F^\times$. In order to use matrices in our
study of gradings we fix the following basis of $\k$:

$$B=(1, e\ot e, x\ot x, x\ot y, y\ot x, y\ot y,
e\ot x, x\ot e, e\ot y, y\ot e)$$
in which the first six elements span the even part while the four last
elements span the odd part. Taking into account that any $t_\lambda$ fixes $e$
and $t_\lambda(x)=\lambda x$, $t_\lambda(y)=\lambda^{-1}y$, the matrix
of $t_{\l,\m}$ relative to $B$ is
$$\hbox{diag}\Big(1,1,\l\m,\frac{\l}{\m},\frac{\m}{\l},\frac{1}{\l\m},
\m,\l,\frac{1}{\m},\frac{1}{\l}\Big).$$

Furthermore, as  conjugated elements produce isomorphic gradings,
we must devote a few lines to the action of the group
$\WW:=\N(\T^2)/\T^2$ on the maximal torus $\T^2$. First of all the
normalizer $\N(\T^2)$ acts on $\T^2$ by conjugation: there is an
action $\N(\T^2)\times\T^2\to \T^2$ such that $(f,t)\mapsto f\circ
t:=ftf^{-1}$, for $f\in\N(\T^2)$ and $t\in\T^2$. Besides the
element $ftf^{-1}$ does not change if we replace $f$ by
$g\in\N(T^2)$ such that $fg^{-1}\in\T^2$. So the previous action
induces an action of $\WW$ on $\T^2$ by conjugation. If
$t,t'\in\T^2$ are in the same orbit under the action of $\WW$ we
shall write $t\sim t'$. Thus we have
\begin{equation}
t_{\l,\m}\sim t_{\l^{-1},\m}\sim t_{\l,\m^{-1}}\sim
t_{\l^{-1},\m^{-1}} \sim t_{\m,\l}
\end{equation}
To prove this, note that $(1,\s)\in\N^2\subset \N(\T^2)$ since
$\N=\T\cup \T\sigma$. But
$$(1,\s)\circ t_{\l,\m}:=(1,\s) t_{\l,\m} (1,\s^{-1})=(1,\s) (t_\l,t_\m) (1,\s^{-1})=$$ $$(t_\l,\s t_\m \s^{-1})=(t_\l,t_{\m^{-1}}\s
 \s^{-1})=t_{\l,\m^{-1}}$$
and similarly $(\s,1)\circ t_{\l,\m}=t_{\l^{-1},\m}$. On the other
hand it is easy to check that $\delta\in\N(T^2)$ satisfies
$\delta\circ t_{\l,\m}=t_{\m,\l}$ so that $t_{\l,\m}\sim
t_{\m,\l}$.
\smallskip

The first step in our study of toral gradings is to look at those
induced by only one toral element $t_{\l,\m}\in\T^2$. It turns out
that this kind of gradings provides most of the cases appearing in
our classification.

\subsubsection{Cyclic gradings}

A {\em cyclic grading} is a toral grading produced by a single
toral element $t_{\l,\m}$. In this case the grading is always
equivalent to a grading by a cyclic group, although not
necessarily the universal group is cyclic (see \cite{G2} for the
concept of universal group.) In order to study the grading induced
by $t_{\l,\m}$ on $J=\k$, which is  the decomposition of $\k$ as a
direct sum of eigenspaces of such toral element $t_{\l,\m}$, we
define the set of eigenvalues
$S:=\{1,\l\m,\frac{\l}{\m},\frac{\m}{\l},\frac{1}{\l\m},
\m,\l,\frac{1}{\m},\frac{1}{\l}\}$ of $t_{\l,\m}$ and consider the
different possibilities for the cardinal $\vert S\vert$. In case
$\vert S\vert=9$ we get the fine toral grading $J=J^1\oplus
J^\l\oplus J^\m\oplus J^{\l\m}\oplus J^{\l\mu^{-1}}\oplus
J^{\m\l^{-1}}\oplus J^{\l^{-1}}\oplus J^{\m^{-1}}\oplus
J^{\l^{-1}\m^{-1}}$ where
\begin{equation}\label{gr1}
\begin{matrix} J^1=\span{1,\ee}, \ J^{\l\m}=\span{\xx},\
J^{\l\m^{-1}}=\span{\xy}, \ J^{\m\l^{-1}}= \span{\yx},\cr
J^{(\l\m)^{-1}}=\span{\yy},
J^\m=\span{\ex},\ J^\l=\span{\xe},\
J^{\m^{-1}}=\span{\ey},\ J^{\l^{-1}}=\span{\ye}.
\end{matrix}\end{equation} This is a $\Z\times\Z$-grading of type $(8,1)$
(that is, eight $1$-dimensional homogeneous components, and one of
dimension $2$), given by $J=\oplus_{(n,m)\in\Z\times\Z}J^{(n,m)}$
with $J^{(n,m)}:=J^{\lambda^n \mu^m}$. Our interest to see this as
a $\Z\times \Z$-grading comes from the fact that this is the
universal grading group, which presents some advantages when is
compared to other groups producing an equivalent grading.  In the
remaining gradings in this section, we will also prefer to use the
superindices to indicate the eigenvalue of the automorphism,
instead of the element of the grading group (of course, both
notations are closely related.)

If $\vert S\vert<9$ we have several possibilities.
\begin{itemize}
\item $1\in\{\l\m,\frac{\l}{\m},\frac{\m}{\l},\frac{1}{\l\m},
\m,\l,\frac{1}{\m},\frac{1}{\l}\}$. Then the grading automorphism
is $t_{\l,\l^{-1}}$, $t_{\l,\l}$, $t_{\l,1}$, or $t_{1,\m}$. But
$t_{\l,\l^{-1}}\sim t_{\l,\l}$ and $t_{\l,1}\sim t_{1,\l}$ so
that, modulo the action of $\WW$, the grading is the one induced
by $t_{\l,\l}$ or $t_{\l,1}$. In case the grading is induced by
$t_{\l,\l}$ we have to distinguish several cases. If the elements
$\l^{-2},\l^{-1},1,\l,\l^2$ are all different we have the grading
 $J=J^{\l^{-2}}\oplus J^{\l^{-1}}\oplus J^1\oplus J^\l\oplus J^{\l^2}$
where
\begin{equation}\label{gr2}
\begin{matrix} J^{\l^{-2}}=\span{\yy},& \quad
J^{\l^{-1}}=\span{\ey,\ye}, &\quad J^1=\span{1,\ee,\xy,\yx},\cr
 J^{\l^{2}}=\span{\xx},& \quad J^{\l}=\span{\ex,\xe}. & &\end{matrix}
\end{equation}
This is a $\Z$-grading of type $(2,2,0,1)$. If there are
coincidences in the set $\{\l^{-2},\l^{-1},$ $1,\l,\l^2\}$ then
necessarily $\lambda$ is a primitive $n$-th root of the unit for
$n=1,2,3,4$. The case $n=1$ corresponds to $\lambda=1$ and thus to
the trivial grading given by $J^1=J$. For $n=2$ we have $\l=-1$ so
that we have a $\Z_2$-grading $J=J^{1}\oplus J^{-1}$ given by
\begin{equation}\label{gr3}
J^1=\span{1,\ee,\xy,\yx,\xx,\yy},\quad
J^{-1}=\span{\ey,\ye,\ex,\xe}.
\end{equation}
This is the $\Z_2$-grading associated to the superalgebra
structure of $J$, of type $(0,0,0,$ $1,0,1)$. For $n=3$, there is
primitive cubic root of $1$, say $\omega$, such that $\l=\omega$.
Then $\lambda^{-1}=\lambda^2$ so that the grading is
$J=J^{\omega^2}\oplus J^1\oplus J^\omega$ where
\begin{equation}\label{gr4}
J^{\omega^2}=\span{\xx,\ey,\ye},\quad  J^1=\span{1,\ee,\xy,\yx},\quad
J^{\omega}=\span{\yy,\ex,\xe}.
\end{equation}
This  is a $\Z_3$-grading of type $(0,0,2,1)$. For $n=4$ we take
$\lambda=i$ the complex unit so that the grading is
$J=J^{-1}\oplus J^{-i}\oplus J^1\oplus J^i$, which is the
$\Z_4$-grading of type $(0,3,0,1)$ given by
\begin{equation}\label{gr5}
\begin{matrix}
J^{-1}=\span{\xx,\yy},& J^{-i}=\span{\ey,\ye},\cr
 J^1=\span{1,\ee,
\xy,\yx}, &  J^i=\span{\ex,\xe}.
\end{matrix}
\end{equation}
In case the grading is induced by $t_{\l,1}$ with $\lambda\ne\pm
1$ we have $J=J^{\l^{-1}}\oplus J^1\oplus J^\l$ which is the
$\Z$-grading of type $(0,0,2,1)$ given by
\begin{equation}\label{gr6}
J^{\lambda^{-1}}=\span{\yx,\yy,\ye}, \ J^1=\span{1,\ee,\ex,\ey},\
J^\l=\span{\xx,\xy,\xe}.
\end{equation}
 For
$\lambda=-1$ we get the coarsening $J=J^{-1}\oplus J^1$ which is a
$\Z_2$-grading of type $(0,0,0,1,0,1)$ given by
\begin{equation}\label{gr7}
J^1=\span{1,\ee,\ex,\ey},\quad
J^{-1}=\span{\xx,\xy,\yx,\yy,\xe,\ye}.
\end{equation}

\item   $1\not\in\{\l\m,\frac{\l}{\m},\frac{\m}{\l},\frac{1}{\l\m},
\m,\l,\frac{1}{\m},\frac{1}{\l}\}$. Now we analyze the possibility
$\l\m\in\{\frac{\l}{\mu},
\frac{\m}{\l},\frac{1}{\l\m},$ $\m,\l,\frac{1}{\m},\frac{1}{\l}\}$.
Modulo the action of the group $\WW$, this gives the following
cases:
\begin{description}

\item{i)} $\l\m=\frac{\l}{\m}$ implying $\m^2=1$. The solution
$\mu=1$ has been considered before. So we have to deal with the
toral element $t_{\l,-1}$. For $\l\ne\pm 1,\pm i$, this induces
the $\Z_2\times\Z$-grading $J=J^1\oplus J^\l\oplus J^{-1}\oplus
J^{-\l}\oplus J^{1/\l}\oplus J^{-1/\l}$ where
\begin{equation}\label{gr8}
\begin{matrix} J^1=\span{1,\ee},\quad J^{-1}=\span{\ex,\ey},\quad
J^\l=\span{\xe},\cr \vbox{\vskip 0.4cm}
J^{-\l}=\span{\xx,\xy},\quad J^{1/\l}=\span{\ye},\quad
J^{-1/\l}=\span{\yy,\yx}, \end{matrix}
\end{equation}
which is of type $(2,4)$. For $\l=1$ the grading is the induced by
$t_{1,-1}\sim t_{-1,1}$, that is, (\ref{gr7}). For $\l=-1$ the
grading is produced  by $t_{-1,-1}$ and this is (\ref{gr3}). For
$\l=i$ we get the grading induced by $t_{i,-1}$, that is:
$J=J^1\oplus J^{-1}\oplus J^i\oplus J^{-1}$. This is the
$\Z_4$-grading of type $(0,2,2)$ such that
\begin{equation}\label{gr9}
\begin{matrix}
J^1=\span{1,\ee},& J^{-1}=\span{\ex,\ey},\cr
 J^i=\span{\xe,\yy,\yx}, &
 J^{-i}=\span{\ye,\xx,\xy}.
\end{matrix}
\end{equation}
For $\l=-i$ we have $t_{-i,-1}$, which is in the orbit of $
t_{i,-1}$ by $\WW$.

\item{ii)}
$\l\m=\frac{1}{\l\m}$. We consider the case $\l=-\frac{1}{\m}$
since the possibility $\lambda\mu=1$ has been studied above. Thus
the grading automorphism is $t_{\l,-\frac{1}{\l}}\sim t_{\l,-
\l}$. If all the elements in the set $\{\pm 1,\pm\l,
\pm\l^{-1},-\l^2,-\l^{-2}\}$ are different, the automorphism
$t_{\l,- \l} $  induces a $\Z_2\times\Z$-grading of type $(6,2)$
given by $J=J^1\oplus J^{-1} \oplus J^{\l}\oplus J^{-\l}\oplus
J^{1/\l}\oplus J^{-1/\l}\oplus J^{-\l^2}\oplus J^{-1/\l^2}$ where
\begin{equation}\label{gr10}
\begin{matrix} J^1=\span{1,\ee},\quad J^{-1}=\span{\xy,\yx},\quad
J^\l=\span{\xe},\quad J^{-\l}=\span{\ex},\cr \vbox{\vskip 0.4cm}
J^{1/\l}=\span{\ye},\quad J^{-1/\l}=\span{\ey}, \quad
J^{-\l^2}=\span{\xx},\quad J^{-1/\l^2}=\span{\yy}.\quad \end{matrix}
\end{equation}
If there are coincidences in the above set the possibilities are:
$\l=\pm 1$, $\l^2=- 1$ or $\l^3=\pm 1$. The first two cases have
been previously studied. The possibility $\l^3=1$ gives
$\l=\omega$ a primitive cubic root of the unit so that this is a
$\Z_6$-grading of type $(2,4)$ induced by $t_{\omega,-\omega}$ and
given by $J=J^1\oplus J^{-1}\oplus J^\omega\oplus
J^{-\omega}\oplus J^{\omega^2}\oplus J^{-\omega^2}$ and
\begin{equation}\label{gr11}
\begin{matrix}
J^1=\span{1,\ee},\quad J^{-1}=\span{\xy,\yx},\quad J^\omega=\span{\xe},\cr
\vbox{\vskip 0.4cm}
\quad J^{-\omega}=\span{\yy,\ex},\quad J^{\omega^2}=\span{\ye},\quad
J^{-\omega^2}=\span{\xx,\ey}.\end{matrix}
\end{equation}
And the possibility $\l^3=-1$   gives $\l=-\omega$, but   $
t_{-\omega,\omega}\sim t_{\omega,-\omega}$, which has just been
  studied.

\item{iii)} $\lambda \mu=\mu$. This would imply $\lambda =1$.

\item{iv)} $\lambda \mu=\lambda$. This would imply $\mu =1$.

\item{v)} $\lambda \mu=\frac{1}{\lambda}$. Then modulo the action
of $\WW$ the grading automorphism is $t_{\lambda, \lambda^2}$. If
all the elements in the set $\{1,\lambda,\lambda^{-1},
\lambda^2,\lambda^{-2},\lambda^3, \lambda^{-3}\}$ are different,
we get  a $\Z$-grading of type $(4,3)$ given by
$$J=J^{\l^{-3}}\oplus J^{\l^{-2}}\oplus J^{\l^{-1}}\oplus J^{1}
\oplus J^{\l}\oplus J^{\l}\oplus J^{\l^{2}}\oplus J^{\l^{3}}$$
where
\begin{equation}\label{gr15}
 \begin{matrix} J^{\l^{-3}}=\span{\yy},\quad J^{\l^{-2}}=\span{\ey},\quad
J^{\l^{-1}}=\span{\xy,\ye},\cr \vbox{\vskip 0.4cm}  \quad
J^{1}=\span{1,\ee},\quad J^{\l}=\span{\yx, \xe},\quad
J^{\l^2}=\span{\ex}, \quad J^{\l^3}=\span{\xx}.\end{matrix}
\end{equation}
If there are coincidences in the above set, the possibilities are:
$\l=\pm 1$, $\l^3= 1$, $\l^4= 1$, $\l^5= 1$ or $\l^6=1$. The first
possibility has been previously considered. The possibility
$\l^3=1$ gives $\l=\omega$ a primitive cubic root of the unit so
that the grading is induced by
$t_{\omega,\omega^2}=t_{\omega,\omega^{-1}} \sim
t_{\omega,\omega}$, previously studied.  The possibility $\l^4=1$
implies $\l \in \pm 1, \pm i$ and so it has also been studied. Let
us consider the case $\l^5= 1$. Hence  $\l=\kappa$ is a primitive
fifth root of the unit. We obtain a $\Z_5$-grading of type $(0,5)$
induced by $t_{\kappa,\kappa^2}$, given by $J=J^{\kappa^4}\oplus
J^{\kappa^3}\oplus J^{\kappa^2} \oplus J^{\kappa} \oplus J^1$
where
\begin{equation}\label{gr17}
\begin{matrix} J^{\kappa^4}=\span{\xy, \ye},\quad
J^{{\kappa^3}}=\span{\xx,\ey},\quad
J^{{\kappa^2}}=\span{\yy,\ex},\cr \vbox{\vskip 0.4cm} \quad
J^{{\kappa}}=\span{\yx,\xe},\quad J^{1}=\span{1,\ee}.\end{matrix}
\end{equation}
Finally, we have to investigate the possibility $\l^6= 1$ which
gives $\l=-\omega$ a primitive sixth root of the unit. We get the
$\Z_6$-grading induced by $t_{-\omega,\omega^2}\sim
t_{-\omega,\frac1{\omega^2}} t_{-\omega, \omega}$, and hence equivalent to (\ref{gr11}).
\end{description}
\end{itemize}
Thus we finish the analysis of the possibility $\l\m\in\{\frac{\l}{\m},
\frac{\m}{\l},\frac{1}{\l\m},\m,\l,\m^{-1},\l^{-1}\}$. Therefore
we must continue analyzing the cases coming from
$\frac{\l}{\m}\in\{ \frac{\m}{\l},\frac{1}{\l\m},\m,\l,\m^{-1},\l^{-1}\}$.
However this analysis reveals no new possibilities for cyclic gradings.
To summarize the results in this subsection we have:

\begin{te}
The nontrivial cyclic gradings on $\k$ are those described in
(\ref{gr1})-(\ref{gr17}).
\end{te}
We list now the essential information of cyclic gradings on the
following table:
\medskip

\vbox{
\hskip 1.5cm
\hbox{
\begin{tabular}{|c|c|c|c|}
\hline Number & Group & type & quasitorus \cr \hline\hline
(\ref{gr1}) & $\Z\times\Z$ & $(8,1)$ & $\T$\cr (\ref{gr2}) & $\Z$
& $(2,2,0,1)$ & $\span{ t_{\l,\l}\colon \lambda\in F^\times}$\cr
(\ref{gr3}) & $\Z_2$ & $(0,0,0,1,0,1)$ & $\span{t_{-1,-1}}$\cr
(\ref{gr4}) & $\Z_3$ & $(0,0,2,1)$ & $\span{t_{\omega,\omega}}$\cr
(\ref{gr5}) & $\Z_4$ & $(0,3,0,1)$ & $\span{t_{i,i}}$\cr
(\ref{gr6}) & $\Z$ & $(0,0,2,1)$ & $\span{t_{\l,1}\colon
\lambda\in F^\times}$\cr 
(\ref{gr7}) & $\Z_2$ & $(0,0,0,1,0,1)$ & $\span{t_{-1,1}}$\cr
(\ref{gr8}) & $\Z\times\Z_2$ & $(2,4)$ & $\span{t_{\l,-1}\colon
\lambda\in F^\times}$\cr
(\ref{gr9}) & $\Z_4$ & $(0,2,2)$ & $\span{t_{i,-1}}$ \cr
(\ref{gr10}) & $\Z\times\Z_2$ & $(6,2)$ & $\span{t_{\l,-\l}\colon
\l\in F^\times}$ \cr
(\ref{gr11}) & $\Z_6$ & $(2,4)$ & $\span{t_{\omega,-\omega}}$ \cr
(\ref{gr15}) & $\Z$ & $(4,3)$ & $\span{t_{\l,\l^2}\colon\l\in F^\times}$ \cr
(\ref{gr17}) & $\Z_5$ & $(0,5)$ & $\span{t_{\kappa,\kappa^2}}$ \cr
\hline
\end{tabular}
\smallskip

} } \noindent where $\omega$, $i$, and $\kappa$ are respectively
$3^{\hbox{\small rd}}$, $4^{\hbox{\small th}}$ and
$5^{\hbox{\small th}}$ primitive roots of the unit. It should be
noted that we have chosen as \textsl{quasitorus} in the table the
maximal quasitorus $Q$ containing $t_{\lambda,\mu}$ but producing
the same grading than the automorphism $t_{\lambda,\mu}$ alone.
For that choice, the group of characters
$\chi(Q)=\rm{Hom}(Q,F^\times)$ is just the universal group.

Note also that the gradings (\ref{gr4}) and (\ref{gr6}), and
(\ref{gr8}) and (\ref{gr11}), are not equivalent in spite of being
of the same type, because their universal groups are not
isomorphic. The gradings (\ref{gr3}) and (\ref{gr7}) are also
nonequivalent, because the automorphisms $t_{-1,-1}$ and
$t_{-1,1}$ are not conjugated in $\aut(\k)$. We can discard also
the equivalence of these gradings by a dimensional argument.

\subsubsection{Noncyclic toral gradings}
To determine the rest of the toral gradings,  we have to continue
by studying the possible refinements of cyclic toral gradings of
the last section. To do that, we consider the grading automorphism
$t_{\lambda, \mu}$ in each case and analyze a refinement to a
grading induced by a set of automorphisms $\{ t_{\lambda, \mu},
t_{\alpha, \beta} \}$. Ruling out the fine grading
 ({\ref{gr1}}), which
 cannot be further refined, the rest of the cyclic gradings do not provide cyclic
refinements except in two cases:

\begin{enumerate}
\item {\it Refinements of ({\ref{gr3}}).} The grading is produced by $t_{-1,
-1}$. An analysis as before of all proper refinements, yields
either cyclic gradings or the $\Z_2 \times \Z_2$-grading induced
by $\{t_{-1,-1}, t_{-1,1}\}$ and given by $J=J^{1,1} \oplus
J^{1,-1} \oplus J^{-1,1}\oplus J^{-1,-1}$ where
\begin{equation}\label{gr19}
\begin{matrix} J^{1,1}=\span{1,\ee},\quad J^{1,-1}=\span{\xx,\xy,
\yx,\yy}, \cr \vbox{\vskip 0.4cm} \quad J^{-1,1}=\span{\ex,\ey},
\quad J^{-1,-1}=\span{\xe,\ye}.\end{matrix}
\end{equation}
This grading is of type $(0,3,0,1)$.  Though its type is the same
as ({\ref{gr5}}), both gradings are not equivalent because the
universal groups are different.

\item {\it Refinements of ({\ref{gr5}}).} The grading is produced
by
 $t_{i, i}$ with $i$ the complex unit. The refinement induced by
 $\{t_{i,i}, t_{-1,1}\}$ is a
$\Z_4 \times \Z_2$-grading given by
$$J=J^{-1,-1} \oplus J^{-i,1} \oplus J^{-i,-1}\oplus J^{1,1} \oplus J^{1,-1}
\oplus J^{i,1} \oplus J^{i,-1}$$ where

\begin{equation}\label{gr23}
\begin{matrix} J^{-1,-1}=\span{\xx,\yy},\quad J^{-i,1}=\span{\ey}, \quad
J^{-i,-1}=\span{\ye}, \quad J^{1,1}=\span{1, \ee} \cr \vbox{\vskip
0.4cm} \quad J^{1,-1}=\span{\xy,\yx}, \quad J^{i,1}=\span{\ex},
\quad J^{i,-1}=\span{\xe}.\end{matrix}
\end{equation}
 This grading is of type (4,3), although not equivalent to the grading ({\ref{gr15}}),
 again because they have nonisomorphic universal groups.
\end{enumerate}

So far, we have only detected two noncyclic gradings produced by
two toral elements which refine the cyclic gradings. To complete
our study of refinements we must describe now the possible
refinements induced by a set of automorphisms $\{t_{-1,-1},
t_{-1,1}, t_{\epsilon,\gamma} \}$ and $\{t_{i,i}, t_{-1,1},
t_{\epsilon,\gamma} \}$. But again a straightforward analysis of
the different possibilities reveals the inexistence of new
gradings. Thus the unique proper refinements (up to equivalences)
of gradings are the ones given in ({\ref{gr19}}) and
({\ref{gr23}}).

\begin{te}
The nontrivial toral  gradings on $\k$ are those described in
(\ref{gr1})-(\ref{gr23}).
\end{te}
The following table contains the relevant information on all toral gradings:
\medskip

\vbox{
\hskip 1.5cm
\hbox{
\begin{tabular}{|c|c|c|c|}
\hline Number & Group & type & quasitorus \cr \hline\hline
(\ref{gr1}) & $\Z\times\Z$ & $(8,1)$ & $\T$\cr (\ref{gr2}) & $\Z$
& $(2,2,0,1)$ & $\span{ t_{\l,\l}\colon \lambda\in F^\times}$\cr
(\ref{gr3}) & $\Z_2$ & $(0,0,0,1,0,1)$ & $\span{t_{-1,-1}}$\cr
(\ref{gr4}) & $\Z_3$ & $(0,0,2,1)$ & $\span{t_{\omega,\omega}}$\cr
(\ref{gr5}) & $\Z_4$ & $(0,3,0,1)$ & $\span{t_{i,i}}$\cr
(\ref{gr6}) & $\Z$ & $(0,0,2,1)$ & $\span{t_{\l,1}\colon
\lambda\in F^\times}$\cr (\ref{gr7}) & $\Z_2$ & $(0,0,0,1,0,1)$ &
$\span{t_{-1,1}}$\cr 
(\ref{gr8}) & $\Z\times\Z_2$ & $(2,4)$ & $\span{t_{\l,-1}\colon
\lambda\in F^\times}$\cr (\ref{gr9}) & $\Z_4$ & $(0,2,2)$ &
$\span{t_{i,-1}}$ \cr (\ref{gr10}) & $\Z\times\Z_2$ & $(6,2)$ &
$\span{t_{\l,-\l}\colon \l\in F^\times}$ \cr (\ref{gr11}) & $\Z_6$
& $(2,4)$ & $\span{t_{\omega,-\omega}}$ \cr (\ref{gr15}) & $\Z$ &
$(4,3)$ & $\span{t_{\l,\l^2}\colon\l\in F^\times}$ \cr
(\ref{gr17}) & $\Z_5$ & $(0,5)$ & $\span{t_{\kappa,\kappa^2}}$ \cr
(\ref{gr19}) & $\Z_2\times\Z_2$ & $(0,3,0,1)$ & $\span{t_{-1,-1},
t_{-1,1}}$\cr (\ref{gr23}) & $\Z_4\times\Z_2$ & $(4,3)$ &
$\span{t_{i,i}, t_{-1,1}}$\cr \hline
\end{tabular}
}
}

\subsection{Fine  gradings}

The fine gradings on $\k$ are induced by maximal abelian subgroups
of diagonalizable automorphisms, MAD-groups from now on,  of
$\aut(\k)$. These are contained in the normalizer of some maximal
torus of $\aut(\k)$. Hence, up to conjugacy, any MAD-group of
$\aut(\k)$ is contained in $\N(\T^2)=\N^2\cup\N^2\delta$, where
$\N$ is the normalizer of $\T$ in $\SL_2(F)$. We refer the reader
to Section II and III for notations. Recall also from previous
sections the following observations: (i) For $(f,g)\in
\SL_2(F)^2$, $\delta (f,g)\delta^{-1} =(g,f)$. (ii) For any
$t_{\lambda} \in \T$, $\sigma t_{\lambda}
\sigma^{-1}=t_{\frac1\lambda}$; and (iii) For any $t \in \T$,
there exists $p\in \SL_2(F)$ such that $p(t \sigma) p^{-1} \in \T$
(which follows from the connectedness of $\SL_2(F)$.)

As examples of MAD-groups of $\aut(\k)$ we have $\T^2$ and

\centerline{${\mathcal M}:= \{(t_{\lambda}, t_{\lambda})\}\times
\{1, \delta\} \subset \N(\T ^2).$}

 \noindent Indeed, if $(f,g)\in\SL_2(F)^2$ commutes with $\mathcal{M}$, it also
does with $\delta$, and $f=g$ according to i). But it also
commutes with $(t,t)\in\T^2$, so $f\in\rm Z_{\rm{SL}_2(F)}(
\T)=\T$, where $Z_G(H)$ denotes the centralizer of $H$ in $G$. Let
us see that, up to conjugacy, these are the only examples.

\begin{te} The MAD-groups of $\aut(\k)$ are, up to conjugacy, $\T^2$
and ${\mathcal M}$.
\end{te}

Proof. Let $A$ be a MAD-group  of $\aut(\k)$, which can be  taken
contained in  $\N(\T^2)=\N^2\cup\N^2\delta$. First suppose that $A
\subset \N^2=\T^2 \cdot \{1,(\sigma,1), (1, \sigma), (\sigma,
\sigma)\}$.
\begin{itemize}
\item If $A \subset \T^2$, then $A=\T^2$.
\item If $A \cap (\T^2 \cdot (\sigma,1))\neq \emptyset$, there exist $t_1, t_2 \in \T$ such that  $(t_1\sigma , t_2
)\in A$. If $t_2\ne\pm1$, $A\subset {\rm Z}_{{\N^2}}(t_1\sigma ,
t_2)=\span{(t_1\sigma,t):t\in\T}$. As there is
  some $p\in\SL_2(F)$ such that $p(
t_1\sigma)p^{-1}=t_3\in\T$, thus
$(p,1)A(p^{-1},1)\subsetneqq\T^2$, what contradicts the maximality
of $A$. If $t_2=\pm1$, then $A\subset {\rm Z}_{{\N^2}}(t_1\sigma ,
\pm1)=\span{(t_1\sigma,t):t\in\T}\cdot\span{(1,\sigma)},$ and
$(p,1)A(p^{-1},1)\subset\span{(t_3,t):t\in\T}\cdot\span{(1,\sigma)}$,
necessarily containing an element of the form $(t_3^n,t_4\sigma)$
with $n\in\mathbb N$ and $t_4\in\T$. Now take $q\in\SL_2(F)$ such
that $q( t_4\sigma)q^{-1}\in\T$ and so
$(p,q)A(p,q)^{-1}\subsetneqq\T^2$, again a contradiction.
\item If  $A \cap (\T^2 \cdot (1,\sigma))\neq
\emptyset$,  there exist $t_1, t_2 \in \T$ such that  $(t_1 , t_2
\sigma)\in A$. As $\delta (t_1 , t_2 \sigma) \delta^{-1}=(t_2
\sigma, t_1)$, we have   $(\delta A \delta^{-1})\cap (\T^2 \cdot
(\sigma,1))\neq \emptyset$, what, taking into account the previous
case, is a contradiction.
\item If $A \cap (\T^2 \cdot
(\sigma,\sigma))\neq \emptyset$, there exist $t_1, t_2 \in \T$
such that  $(t_1 \sigma, t_2 \sigma)\in A$. As any other element
of $A$ commutes with it, it is of some of the following types:
$(\pm t_1\sigma,\pm t_2\sigma)$, $(\pm 1,\pm t_2\sigma)$, $(\pm
t_1\sigma,\pm 1)$ or $(\pm 1,\pm 1)$. Thus taking $p,q\in
\SL_2(F)$ such that $p(t_1\sigma)p^{-1}$,
$q(t_2\sigma)q^{-1}\in\T$, we find that
$(p,q)A(p,q)^{-1}\subsetneqq\T^2$ and again the contradiction
appears.
\end{itemize}
We have shown  that in the case $A \subset \N^2$, necessarily
$A=\T^2$. Suppose next that $A\nsubseteq \N^2$, so that $A \cap
\N^2 \delta \neq \emptyset$. Then there exist $t_1, t_2 \in \T$
such that either $(t_1 , t_2 ) \delta\in A$, $(t_1 \sigma , t_2 )
\delta\in A$, $(t_1 , t_2 \sigma) \delta\in A$ or $(t_1 \sigma ,
t_2 \sigma) \delta\in A$. Observe that the third possibility can
be reduced to the second one by conjugating with $\delta$, and
that the forth possibility can be reduced to the first one by
conjugating with $(\sigma,1)$. Besides the second possibility can
be reduced to the first one by conjugating with
$(pt_2^{-1},p)\delta$, for $p\in\SL_2(F)$ such that
$p(t_1t_2^{-1}\sigma)p^{-1}\in\T$.

So we will consider the first possibility. We can suppose that
certain  $(1,t_1)\delta \in A$ by conjugating $(t_1,t_2)\delta\in
A$  with $(t_1^{-1},1)$. Take $s\in \T$ such that $s^2=t_1$.
Because of the abelian character of $A$, it is contained in  $
{\rm Z}_{\N(\T^2)}((1,t_1) \delta)$. In case $t_1 \neq \pm1$, we
have
 ${\rm Z}_{{\N(T^2)}}((1,t_1)
\delta)=\{(t,t):t\in \T\} \cdot \{1, (1,t_1) \delta\}$, and hence
   $(s,1)A(s,1)^{-1} \subset \{(t,t):t\in \T\}
\cdot \{1, (s,s) \delta\}=\mathcal{M}$, so that the equality holds
by maximality. In the case $t_1 = \pm 1$, one has ${\rm Z}_{
{\N(\T^2)}}((1,\pm 1) \delta) =\{(t,t):t\in \T\} \cdot\span{
(1,\pm 1) \delta}\cdot \span{(\sigma, \sigma)}:=\widetilde{A}.$ If
$A\subset\{(t,t):t\in \T\} \cdot\span{ (1,\pm 1) \delta}$, as
before $(s,1)A(s,1)^{-1}\subset\mathcal{M}$ and is equal to
$\mathcal{M}$ by maximality. Otherwise, there exists $t_2 \in \T$
such that $(t_2\sigma, t_2 \sigma) \in A$. Then $A \subset {\rm
Z}_{\widetilde{A}}(t_2 \sigma, t_2
\sigma)=\span{(1,\pm1)\delta}\cdot\span{(t_2 \sigma, t_2
\sigma)}$. For $t_1=1$, take  $p \in \SL_2(F)$ such that $p(t_2
\sigma) p^{-1} \in \T$ and so  $(p,p)A(p,p)^{-1} \subsetneqq
{\mathcal {M}}$, what contradicts the maximality of $A$. For
$t_1=-1$, take $p \in \SL_2(F)$ such that $p(\sigma t_2 ) p^{-1}
\in \T$ and check that
$(p\sigma,pt_2^{-1})A(p\sigma,pt_2^{-1})^{-1} \subsetneqq
{\mathcal {M}}$, again a contradiction.

\begin{co}
The fine gradings on $\k$ are, up to equivalence:
\begin{itemize}
\item[{\rm(1)}] The toral $\Z\times\Z$-grading,
$J^{(0,0)}=\span{1,\ee}, J^{(1,1)}=\span{\xx},
J^{(1,-1)}=\span{\xy}, J^{(-1,1)}= \span{\yx},
 J^{(-1,-1)}=\span{\yy}, \\
  J^{(0,1)}=\span{\ex},
J^{(1,0)}=\span{\xe}, J^{(0,-1)}=\span{\ey},
J^{(-1,0)}=\span{\ye}.$ This is of type $(8,1)$.

\item[{\rm(2)}]  The non-toral $\Z \times \Z_2$-grading, $
J^{(0,{0})}=\span{1,\ \ee,\ \xy-\yx}, J^{(-1,{1})} =\span{\ey-\ye},
J^{(2,{1})}=\span{\xx}, J^{(-2,{1})}=\span{\yy},
J^{(1,{1})}=\span{\ex-\xe}, J^{(0,{1})} =\span{\xy+\yx},
J^{(1,{0})}=\span{\ex+\xe}, J^{(-1,{0})}=\span{\ey+\ye}.$ This is
of type $(7,0,1)$.
\end{itemize}
\end{co}

\subsection{Nontoral gradings}

Recall that the quasitorus $\mathcal{M}=\langle t_{\lambda,
\lambda},\delta:\lambda\in F^\times\rangle\cong
F^\times\times\Z_2$ induces a $\Z\times\Z_2$-grading  of type
$(7,0,1)$
 with homogeneous components
 \begin{equation}\label{lafina}
\begin{matrix} J^{(1,1)}=\span{1, \ee, \xy-\yx}, &
J^{(\lambda^2,-1)} =\span{\xx}, \cr
J^{(\frac1{\lambda^2},-1)}=\span{\yy}, &
J^{(\lambda,1)}=\span{\ex+\xe}, \cr
J^{(\frac1\lambda,1)}=\span{\ey+\ye}, & J^{(1,-1)}
=\span{\xy+\yx}, \cr J^{(\lambda,-1)}=\span{\ex-\xe}, &
J^{(\frac1\lambda,-1)}=\span{\ey-\ye},\end{matrix} \end{equation}
 where the coordinates of the superindex indicate now the eigenvalues of the actions of $t_{\lambda, \lambda}$ and $\delta$ respectively.
 Note that any nontoral grading is produced by a subquasitorus $Q\subset \mathcal{M}$ such that $Q\cap\{\delta t_{\lambda, \lambda} :\lambda\in F^\times\}\ne\emptyset$ (otherwise $Q\subset \mathcal{T}':=\{t_{\lambda, \lambda} :\lambda\in F^\times\}$ would be toral.)
Note also that, for any $\delta t_{\beta,\beta}\in Q$,
 $$
 Q=(Q\cap\mathcal{T}')\cdot\langle \delta t_{\beta,\beta}\rangle.
 $$
First suppose that $Q=\langle \delta t_{\beta,\beta}\rangle$.
Recall that $\delta t_{\beta,\beta}$ acts in the homogeneous
components of (\ref{lafina}) with eigenvalues $\{1,
-\beta^2,-\frac1{\beta^2},\beta,\frac1\beta,-1,-\beta,-\frac1\beta\}$
respectively. Hence it produces a grading equivalent to
(\ref{lafina}) if those numbers are all different, that is, if
$\beta^4\ne1$ and $\beta^6\ne1$. Let us see which is the induced
grading for any value of $\beta=\pm1,\pm i,\pm\omega$, for
$\omega$ a cubic root of unit, taking into account that
$t_{1,-1}(\delta t_{\beta,\beta})(t_{1,-1})^{-1}=\delta
t_{-\beta,-\beta}$, and so $\langle \delta t_{\beta,\beta}\rangle$
and $\langle \delta t_{-\beta,-\beta}\rangle$ produce equivalent
gradings.
\begin{itemize}
\item If $\beta=1$, we get the $\Z_2$-grading produced by $\delta$, given by
\begin{equation}\label{z2}\begin{matrix} J^{1}=\span{1,\ \ee,\ \xy-\yx, \
\ex+\xe,\ \ey+\ye}, \quad
 \cr \vbox{\vskip 0.4cm} \quad
J^{-1}=\span{\ex-\xe,\ \ey-\ye,\ \xy+\yx,\ \xx,\ \yy}, \end{matrix}
\end{equation} of type
(0,0,0,0,2).
\item If $\beta=i$, we get the
$\Z_4$-grading  given by
\begin{equation}\label{z4}\begin{matrix} J^{1}=\span{1,\ \ee,\ \xy-\yx,\
\xx, \yy}, \quad J^{i}=\span{\xe+\ex,\ \ye-\ey},
 \cr \vbox{\vskip 0.4cm} \quad
J^{-1}=\span{\xy+\yx}, \quad J^{-i}=\span{\ye+\ey,\ \xe-\ex},
\end{matrix}
\end{equation} of type (1,2,0,0,1).
\item If $\beta=\omega$, we get  the
$\Z_6$-grading   given by
\begin{equation}\label{z6}
\begin{matrix} J^{1}=\span{1,\ \ee,\ \xy-\yx}, & J^{-1}=\span{
\xy+\yx},\cr 
J^{-\omega^2}=\span{\xx,\ \ye-\ey}, & 
J^{ -\omega}=\span{\yy,\ \xe-\ex}, \cr J^{\omega}=\span{\xe+\ex}, & 
J^{\omega^2}=\span{  \ey+\ye}, \end{matrix} \end{equation}

of type (3,2,1).
\end{itemize}
In case that  $Q\ne\langle \delta t_{\beta,\beta}\rangle$,
since $t_{\beta^2,\beta^2}=(\delta t_{\beta,\beta})^2\in
Q\cap\mathcal{T}'$, we get $\langle
t_{\beta^2,\beta^2}\rangle\subsetneqq
 Q\cap\mathcal{T}'$. If there
is some $t_{\lambda,\lambda}\in Q\cap\mathcal{T}'$ with
$\lambda^4$ and $\lambda^3\ne1$, then the grading induced by $Q$
is not a coarsening of (\ref{lafina}), because the pair
$(t_{\lambda,\lambda},\delta t_{\beta,\beta})$ acts with
eigenvalues
$$\{(1,1), (\lambda^2,-\beta^2), (\lambda^{-2},- {\beta^{-2}}),
(\lambda,\beta), (\lambda^{-1},\beta^{-1}),
(1,-1),(\lambda,-\beta), (\lambda^{-1}, -\beta^{-1}) \},$$ which
are all different (in fact the first coordinates  are different,
up to the pairs $(1,1)$ and $(1,-1)$). Hence we can assume that
every element in $Q\cap\mathcal{T}'$ has order divisor of either 3
or 4. As $Q\cap\mathcal{T}'$ is a subgroup of $\mathcal{T}'\cong
F^\times$, we conclude that $Q\cap\mathcal{T}'$ equals  either
$\langle t_{-1,-1}\rangle$, or $\langle t_{i,i}\rangle$, or
$\langle t_{\omega,\omega}\rangle$.
\begin{itemize}
\item If $Q\cap\mathcal{T}'=\langle t_{-1,-1}\rangle$,
then $\beta^2=1$ (since $\beta^2$ does not generate $\{1,-1\}$)
and $Q=Q\cap\mathcal{T}'\cdot\langle\delta
t_{\beta,\beta}\rangle=\langle
t_{-1,-1},\delta\rangle\cong\Z_2\times\Z_2$, which induces a
grading of type $(0,2,2)$:
\begin{equation}\label{Z22}
\begin{matrix} J^{(1,1)}=\span{1,\ \ee,\ \xy-\yx} & J^{(1,-1)}
=\span{\xy+\yx,\ \xx,\ \yy}, \cr J^{(-1,1)}=\span{\ex+\xe,\
\ey+\ye}, &
J^{(-1,-1)}=\span{\ex-\xe, \ \ey-\ye}.\end{matrix}
\end{equation}

\item If $Q\cap\mathcal{T}'=\langle t_{i,i}\rangle$,  then $\beta^2=\pm1$
($\beta^2$ does not generate $\{1,-1,i,-i\}$) so that $\beta\in\{1,-1,i,-i\}$,
$\delta\in Q$ and $Q=\langle t_{i,i},\delta\rangle\cong\Z_4\times\Z_2$, which
induces
\begin{equation}\label{Z42}
\begin{matrix} J^{(1,1)}=\span{1,\ \ee,\ \xy-\yx},& J^{(i,1)}
=\span{\ex+\xe}, \cr J^{(-1,1)}=\span{\xx, \yy}, &
J^{(-i,1)}=\span{\ey+\ye}, \cr  J^{(1,-1)} =\span{\xy+\yx}, &
J^{(i,-1)}=\span{\ex-\xe}, \cr
J^{(-i,-1)}=\span{\ey-\ye},\end{matrix}
\end{equation}
a grading of type $(5,1,1)$.
\item If $Q\cap\mathcal{T}'=\langle t_{\omega,\omega}\rangle$,
then $\beta^2=1$ (the only non-generator in $\Z_3$) so that
$\beta=\pm1$. Both cases are conjugated, since $t_{-1,1}\delta
(t_{-1,1})^{-1}=\delta t_{-1,-1}$. Thus we can assume that
$\delta\in Q$ and $Q=\langle
t_{\omega,\omega},\delta\rangle=\langle\delta
t_{\omega,\omega}\rangle$, so that  we are in the first case.
\end{itemize}

Summarizing the results in these sections, we  claim

\begin{te} Up to equivalence, the nontrivial $G$-gradings on $J=K_{10}$ (compatible with the superalgebra structure) are the
following:

\begin{itemize}

\item[{\rm(1)}] $G=\Z\times\Z$,    $ J^{(0,0)}=\span{1,\ee},
J^{(1,1)}=\span{\xx}, J^{(1,-1)}=\span{\xy},
J^{(-1,1)}=\span{\yx},
 J^{(-1,-1)}=\span{\yy}, \\
  J^{(0,1)}=\span{\ex},
J^{(1,0)}=\span{\xe}, J^{(0,-1)}=\span{\ey},
J^{(-1,0)}=\span{\ye}. $

\item[{\rm(2)}] $G=\Z$, $ J^{{-2}}=\span{\yy},
J^{{-1}}=\span{\ey,\ye}, J^0= \span{1,\ee,\xy,\yx},
 J^{{2}}=\span{\xx}, J^{{1}}=\span{\ex,\xe}.$

 \item[{\rm(3)}] $G=\Z_2$, $J^0=\span{1,\ee,\xy,\yx,\xx,\yy}, J^{1}=\span{\ey,\ye,\ex,\xe}$.

\item[{\rm(4)}] $G=\Z_3$, $J^{2}=\span{\xx,\ey,\ye},
J^0=\span{1,\ee,\xy,\yx}, J^{1}=\span{\yy,\ex,\xe}.$

\item[{\rm(5)}] $G=\Z_4$, $ J^{2}=\span{\xx,\yy},
J^{3}=\span{\ey,\ye}, J^0=\span{1,\ee, \xy,\yx},
J^1=\span{\ex,\xe}.$

\item[{\rm(6)}] $G=\Z$, $J^{{-1}}=\span{\yx,\yy,\ye},
J^0=\span{1,\ee,\ex,\ey}, J^1=\span{\xx,\xy,\xe}.$

\item[{\rm(7)}] $G=\Z_2$,  $J^0=\span{1,\ee,\ex,\ey},
J^{1}=\span{\xx,\xy,\yx,\yy,\xe,\ye}.$

\item[{\rm(8)}] $G=\Z_2\times\Z$, $J^{(0,0)}=\span{1,\ee},
J^{(1,0)}=\span{\ex,\ey}, J^{(0,1)}=\span{\xe},
J^{(1,1)}=\span{\xx,\xy}, J^{(0,-1)}=\span{\ye},
J^{(1,-1)}=\span{\yy,\yx}. $

\item[{\rm(9)}] $G=\Z_4$, $J^0=\span{1,\ee}, J^{2}=\span{\ex,\ey},
J^1=\span{\xe,\yy,\yx}, J^{3}=\span{\ye,\xx,\xy}.$

\item[{\rm(10)}] $G=\Z_2\times\Z$, $ J^{(0,0)}=\span{1,\ee},
J^{(1,0)}=\span{\xy,\yx}, J^{(0,1)}=\span{\xe},
J^{(1,1)}=\span{\ex}, J^{(0,-1)}=\span{\ye},
J^{(1,-1)}=\span{\ey}, J^{(1,2)}=\span{\xx},
J^{(1,-2)}=\span{\yy}.$

\item[{\rm(11)}] $G=\Z_6$, $J^0=\span{1,\ee},
J^{3}=\span{\xy,\yx}, J^2=\span{\xe}, J^{5}=\span{\yy,\ex},
J^{4}=\span{\ye}, J^{1}=\span{\xx,\ey}.$

\item[{\rm(12)}] $G=\Z$, $ J^{{-3}}=\span{\yy},
J^{{-2}}=\span{\ey}, J^{{-1}}=\span{\xy,\ye}, J^{0}=\span{1,\ee},
J^{1}=\span{\yx, \xe}, J^{2}=\span{\ex}, J^{3}=\span{\xx}.$

\item[{\rm(13)}] $G=\Z_5$, $J^{4}=\span{\xy, \ye},
J^{{3}}=\span{\xx,\ey}, J^{{2}}=\span{\yy,\ex},
J^{{1}}=\span{\yx,\xe}, J^{0}=\span{1,\ee}.$

\item[{\rm(14)}] $G=\Z_2 \times \Z_2$, $ J^{(0,0)}=\span{1,\ee},
J^{(0,1)}=\span{\xx,\xy, \yx,\yy}, J^{(1,0)}=\span{\ex,\ey},
J^{(1,1)}=\span{\xe,\ye}.$

\item[{\rm(15)}] $G=\Z_4 \times \Z_2$, $J^{(2,1)}=\span{\xx,\yy},
J^{(3,0)}=\span{\ey}, J^{(3,1)}=\span{\ye},  J^{(0,0)}=\span{1,
\ee}, J^{(0,1)}=\span{\xy,\yx},  J^{(1,0)}=\span{\ex},
J^{(1,1)}=\span{\xe} .$

\item[{\rm(16)}]  $G=\Z \times \Z_2$, $
J^{(0,{0})}=\span{1,\ \ee,\ \xy-\yx},\ J^{(-1,{1})} =\span{\ey-\ye},\
J^{(2,{1})}=\span{\xx}, J^{(-2,{1})}=\span{\yy},\
J^{(1,{1})}=\span{\ex-\xe}, J^{(0,{1})} =\span{\xy+\yx},\
J^{(1,{0})}=\span{\ex+\xe}, J^{(-1,{0})}=\span{\ey+\ye}.$

\item[{\rm(17)}] $G=\Z_4 \times \Z_2$,
$ J^{(0,0)}=\span{1,\ \ee,\ \xy-\yx},\ J^{(1,0)} =\span{\ex+\xe},
J^{(2,0)}=\span{\xx, \yy},\ J^{(3,0)}=\span{\ey+\ye},\  J^{(0,1)}
=\span{\xy+\yx}, J^{(1,1)}=\span{\ex-\xe},\
J^{(3,1)}=\span{\ey-\ye}.$

\item[{\rm(18)}] $G=\Z_2 \times \Z_2$,
$J^{(0,0)}=\span{1,\ \ee,\ \xy-\yx}$, $J^{(0,1)} =\span{\xy+\yx,\
\xx,\ \yy},$ $J^{(1,0)}=\span{\ex+\xe,\ \ey+\ye},$
$J^{(1,1)}=\span{\ex-\xe, \ \ey-\ye}.$

\item[{\rm(19)}] $G=\Z_2$,
$J^{0}=\span{1,\ \ee,\ \xy-\yx, \ \ex+\xe,\ \ey+\ye}$, $ J^{1}=\span{\ex-\xe,\
\ey-\ye,\ \xy+\yx,\ \xx,\ \yy}$.

\item[{\rm(20)}] $G=\Z_4$,
$J^{0}=\span{1,\ \ee,\ \xy-\yx,\ \xx, \yy},$
$J^{1}=\span{\xe+\ex,\ \ye-\ey}$, $J^{2}=\span{\xy+\yx}$,
$J^{3}=\span{\ye+\ey,\ \xe-\ex}$.

\item[{\rm(21)}] $G=\Z_6$,
$J^{0}=\span{1,\ \ee,\ \xy-\yx}$, $J^{1}=\span{\xx,\ \ye-\ey}$,
$J^{2}=\span{\xe+\ex}$, $J^{3}=\span{ \xy+\yx}$, $J^{4}=\span{
\ye+\ey}$,
 $J^{5}=\span{\yy,\
\xe-\ex}$.

The unique fine gradings are {\rm (1)} and {\rm (16)}. The
gradings {\rm (1)}-{\rm (15)} are the toral ones.

\end{itemize}

\end{te}

\end{document}